\documentclass[leqno,12pt]{article}

\usepackage{amsmath,amsthm,amscd,amssymb,graphicx,indentfirst}
\usepackage{amsthm}

\theoremstyle{plain}

\newtheorem*{Theorem}{\indent Theorem}
\newtheorem*{lemma}{\indent Lemma}
\newtheorem*{Cor}{\indent Corollary}

\theoremstyle{definition}
\newtheorem*{acknowledgement}{\indent Acknowledgement}

\makeatletter

\def\address#1#2{\begingroup
\noindent\parbox[t]{7.8cm}{%
\small{\scshape\ignorespaces#1}\par\vskip1ex
\noindent\small{\itshape E-mail address}%
\/: #2\par\vskip4ex}\hfill%
\endgroup}%

\makeatother

\title{Global monodromy modulo 5 of quintic-mirror family}

\author{Kennichiro Shirakawa}
\date{}

\begin{document}

\maketitle


\begin{abstract}
The quintic-mirror family is a well-known one-parameter family of Calabi-Yau threefolds.  A complete description of the global monodromy group of this family is not yet 
known. In this paper, we give a presentation of the global monodromy group in the general linear group of degree 4 over the ring of integers modulo 5. 
\end{abstract}

\section{Introduction}

The quintic-mirror family $(W_\lambda)_{\lambda\in\mathbb{P}^1}\to\mathbb{P}^1$ is a family, whose restriction $f:(W_\lambda)_{\lambda\in U}\to U$ on 
$U:=\mathbb{P}^1-\{0,1,\infty\}$ is a smooth projective family of Calabi-Yau manifolds. Fix $b\in U$ and let $\langle \ , \ \rangle$ be the 
anti-symmetric bilinear form on $H^3(W_b,\mathbb{Z})$ defined by the cup product. The global monodromy group $\Gamma$ is the image of the representation 
$\pi_1(U,b) \to \text{Aut}(H^3(W_b,\mathbb{Z}),\langle \ , \ \rangle)$ corresponding to the local system 
$R^3f_*\mathbb{Z}$ with the fiber $H^3(W_b,\mathbb{Z})$ over $b$. When we take a symplectic basis, we can identify $\text{Aut}(H^3(W_b,\mathbb{Z}),\langle \ , \ \rangle)$ 
with $\text{Sp}(4,\mathbb{Z})$.

In this paper, we are concerned with a description of $\Gamma$. 
Matrix presentations of the generators of $\Gamma$ are well studied and it is also known that $\Gamma$ is Zariski dense in $\text{Sp}(4,\mathbb{Z})$ (e.g. [1], [3]). 
However, it is not known whether the index of $\Gamma$ in $\text{Sp}(4,\mathbb{Z})$ is finite or not (e.g. [2]). A direct approach for this problem is to describe 
$\Gamma$ explicitly. 
In the main theorem of this paper, we give a presentation of $\Gamma$ in $\text{GL}(4,\mathbb{Z}/5\mathbb{Z})$, which is a small attempt toward a description of $\Gamma$. 

On the other hand, Chen, Yang and Yui find a congruence subgroup $\Gamma(5,5)$ of $\text{Sp}(4,\mathbb{Z})$ of finite index, which contains $\Gamma$ in [2]. 
Combining their result and our main theorem, we can construct a smaller congruence subgroup $\tilde{\Gamma}(5,5)$ of $\text{Sp}(4,\mathbb{Z})$ of finite 
index, which contains $\Gamma$. However this result is merely the fact that $\tilde{\Gamma}(5,5)$ contains $\Gamma$. 
After all, the index of $\Gamma$ in $\text{Sp}(4,\mathbb{Z})$ is still unknown.

\section{The quintic-mirror family}

The quintic-mirror family was constructed by Greene and Plesser. We review the construction of the quintic-mirror family after [4].

Let $\psi \in \mathbb{P}^1= \mathbb{C} \cup \{\infty \}$, and let
\[ Q_\psi=\{ x \in \mathbb{P}^4 \ | \ x_1^5+x_2^5+x_3^5+x_4^5+x_5^5-5\psi x_1x_2x_3x_4x_5=0 \}. \]

\smallskip

\noindent A finite group $G$, which is abstractly isomorphic to $(\mathbb{Z}/5\mathbb{Z})^3$, acts on $Q_\psi$ as follows.
\begin{align*}
& \mu_5:\text{the multiplicative group of} \ \text{the} \ 5 \text{-th root of} \ 1 \in \mathbb{C}, \\
& \widetilde{G}=(\mu_5)^5/\{ (\alpha_1,\cdots,\alpha_5) \in (\mu_5)^5 \ | \ \alpha_1=\cdots=\alpha_5 \}, \\
& G=\{ (\alpha_1,\cdots,\alpha_5) \in \widetilde{G} \ | \ \alpha_1 \cdots \alpha_5=1 \}, \\
& G \times Q_\psi \to Q_\psi \ , \ 
( (\alpha_1,\cdots,\alpha_5),(x_1,\cdots,x_5) ) \mapsto ( \alpha_1 x_1,\cdots,\alpha_5 x_5 ).
\end{align*}

\smallskip

When we divide the hypersurface $Q_\psi$ by $G$, canonical singularities appear. 
For $\psi \in \mathbb{C} \subset \mathbb{P}^1,$ it is known that there is a simultaneous minimal desingularization of 
 these singularities, and we have the one-parameter family $(W_\psi)_{\psi \in \mathbb{P}^1}$ whose fibres are listed as follows:

\smallskip
\smallskip
\noindent $\bullet$ When $\psi$ belongs to $\mu_5 \subset \mathbb{C} \subset \mathbb{P}^1$, $W_\psi$ has one ordinary double point.

\smallskip
\smallskip
\noindent $\bullet$ $W_\infty$ is a normal crossing divisor in the total space.

\smallskip
\smallskip
\noindent $\bullet$ The other fibres of $(W_\psi)_{\psi\in \mathbb{P}^1}$ are smooth with Hodge numbers $h^{p,q}=1$ for $p+q=3, \ p,q \ge 0.$

By the action of

\smallskip

\ \ \ \ \ $\alpha \in \mu_5, \ (x_1,\cdots,x_5) \ \mapsto (x_1,\cdots,x_4,\alpha^{-1} x_5),$

\smallskip

\noindent we have the isomorphism from the fibre over $\psi$ to the fibre over $\alpha \psi.$
Let $\lambda=\psi^5$ and let
\[ \begin{array}{ccccc}
		(W_\lambda)_{\lambda\in \mathbb{P}^1} & \resizebox{0.1\hsize}{\height}{=} & ((W_\psi)_{\psi\in \mathbb{P}^1})/\mu_5 \\
		\\
		\big\downarrow & & \big\downarrow \\
		\\
		(\lambda \text{-plane}) & \resizebox{0.1\hsize}{\height}{=} & (\psi \text{-plane})/\mu_5.
			\end{array} \]

\smallskip
\smallskip
\noindent This family $(W_\lambda)_{\lambda\in \mathbb{P}^1}$ is the so-called quintic-mirror family. 
(For more details of the above construction, see e.g. [4], [5].)

\section{Monodromy}

Let $b \in \mathbb{P}^1 - \{0,1,\infty\}$ on the $\lambda$-plane. 
In [1], Candelas, de la Ossa, Green and Parks constructed a symplectic basis $\{A^1,A^2,B_1,B_2\}$ of 
$H_3(W_b,\mathbb{Z})$ and calculated the monodromies around $\lambda=0,1,\infty$ on the period integrals of a holomorphic 3-form on this basis. 
By the relation in [5, Appendix C] between the symplectic basis $\{ \beta^1,\beta^2,\alpha_1,\alpha_2 \}$ of $H^3(W_b,\mathbb{Z})$, 
which is defined to be the dual basis of $\{B_1,B_2,A^1,A^2\}$, and the period integrals, we have the matrix representations of the local monodromies for the basis 
$\{ \beta^1,\beta^2,\alpha_1,\alpha_2 \}$. We recall their results.

Matrix representations $A,T,T_\infty$ of local monodromies around $\lambda=0,1,\infty$ for the basis $\{ \beta^1,\beta^2,\alpha_1,\alpha_2 \}$ are as follows:
\[ A=\begin{pmatrix}
									11 & 8 & -5 & 0 \\
									5 & -4 & -3 & 1 \\
									20 & 15 & -9 & 0 \\
									5 & -5 & -3 & 1
										\end{pmatrix}, \ 
	T=\begin{pmatrix}
											1 & 0 & 0 & 0 \\
											0 & 1 & 0 & -1 \\
											0 & 0 & 1 & 0 \\
											0 & 0 & 0 & 1
												\end{pmatrix}, \ 
T_\infty = \begin{pmatrix}
				-9 & -3 & 5 & 0 \\
				0 & 1 & 0 & 0 \\
				-20 & -5 & 11 & 0 \\
				-15 & 5 & 8 & 1
					\end{pmatrix}. \]
In particular, the above $A$ and $T$ are the inverse matrices of the matrices $A$ and $T$ in the lists of [1] respectively.

Let $\langle \ , \ \rangle$ be the anti-symmetric bilinear form on $H^3(W_b,\mathbb{Z})$ defined by the cup product. 
The global monodromy $\Gamma$ is $\text{Im}(\pi_1(\mathbb{P}^1 - \{0,1,\infty\}) \to \text{Aut}(H^3(W_b,\mathbb{Z}),\langle \ , \ \rangle )$.
When we take $\{ \beta^1,\beta^2,\alpha_1,\alpha_2 \}$ as the basis of $H^3(W_b,\mathbb{Z})$,
$\text{Aut}(H^3(W_b,\mathbb{Z}),\langle \ , \ \rangle )$ is identified with $\text{Sp}(4,\mathbb{Z})$, and $\Gamma$ is the subgroup of 
$\text{Sp}(4,\mathbb{Z})$ which is generated by $A$ and $T$.

We can partially normalize A and T simultaneously as follows.

\begin{lemma}
There exists $P \in \textup{GL}(4,\mathbb{Q})$ such that
\[ P^{-1}A^{-1}P=\begin{pmatrix}
				1 & 1 & 0 & 0 \\
				0 & 1 & 1 & -1 \\
				0 & 0 & 1 & -1 \\
				5 & 5 & 5 & -4
				 \end{pmatrix}, \ 
	P^{-1}T^{-1}P=\begin{pmatrix}
				1 & 0 & 0 & 0 \\
				0 & 1 & 0 & 0 \\
				0 & 0 & 1 & 1 \\
				0 & 0 & 0 & 1
					\end{pmatrix}. \]

\end{lemma}

\begin{proof}
We take $P=\begin{pmatrix}
						5 & -3 & 0 & 0 \\
						0 & 0 & 1 & 0 \\
						10 & -5 & 0 & 0 \\
						0 & 0 & 0 & 1
					\end{pmatrix}$.
The assertion follows.
\end{proof}

\section{Main result}

Let $\Gamma'=\{ P^{-1}XP \in \textup{GL}(4,\mathbb{Z}) \ | \ X \in \Gamma \}$, and let $\rho:\text{GL}(4,\mathbb{Z}) \to \text{GL}(4,\mathbb{Z}/5\mathbb{Z})$ 
be the natural projection. Define $\tilde{\Gamma}=\rho(\Gamma')$. We study $\tilde{\Gamma}$.

Let $\tilde{A}=\rho(P^{-1}A^{-1}P), \ \tilde{T}=\rho(P^{-1}T^{-1}P) \in \text{GL}(4,\mathbb{Z}/5\mathbb{Z})$. By a simple calculation, we obtain
\[ \tilde{A}^{n}=\begin{pmatrix}
						1 & n & 3n(n+4) & n(n+1)(4n+1) \\
						0 & 1 & n & 2n(n+1) \\
						0 & 0 & 1 & 4n \\
						0 & 0 & 0 & 1
					\end{pmatrix} \in \text{GL}(4,\mathbb{Z}/5\mathbb{Z}). \]

Let $\hat{\Gamma}$ be
\[ \left \{ \begin{array}{c|c} 
				\begin{pmatrix}
									1 & n & 3n^2+2n & a \\
									0 & 1 & n & b \\
									0 & 0 & 1 & c \\
									0 & 0 & 0 & 1 
										\end{pmatrix} 
											\in \text{GL}(4,\mathbb{Z}/5\mathbb{Z}) \ & n,a,b,c \in \mathbb{Z}/5\mathbb{Z}
															\end{array} \right \}. \]
$\hat{\Gamma}$ is a subgroup of $\text{GL}(4,\mathbb{Z}/5\mathbb{Z})$ which contains $\tilde{A}$ and $\tilde{T}$. The following Theorem and Corollary are the 
main results of this paper.

\begin{Theorem}
$\tilde{\Gamma}=\hat{\Gamma}$.
\end{Theorem}

\begin{proof}
$\tilde{\Gamma}\subset\hat{\Gamma}$ follows from what we just mentioned. So we shall prove the inverse inclusion.

From the presentations of elements of $\hat{\Gamma}$, we see that $\hat{\Gamma}$ is generated by
\[ \ \tilde{A}, \ \tilde{T}, \ E_1=\begin{pmatrix}
									1 & 0 & 0 & 1 \\
									0 & 1 & 0 & 0 \\
									0 & 0 & 1 & 0 \\
									0 & 0 & 0 & 1
										\end{pmatrix} \ \text{and} \ 
	E_2=\begin{pmatrix}
											1 & 0 & 0 & 0 \\
											0 & 1 & 0 & 1 \\
											0 & 0 & 1 & 0 \\
											0 & 0 & 0 & 1
												\end{pmatrix}. \]
Therefore, it is enough to show $E_1$ and $E_2$ belong to $\tilde{\Gamma}$. In fact, we have
\[ E_2=\tilde{A}\tilde{T}\tilde{A}^4\tilde{T}^4, \ E_1=(E_2^2\tilde{A}^2\tilde{T}^4\tilde{A}^3\tilde{T})^4. \]
Hence $E_1, \ E_2 \in \tilde{\Gamma}$.
\end{proof}

\begin{Cor}
Let $X \in \Gamma$. Then the eigenpolynomial of $X$ is
\[ x^4+(5m+1)x^3+(5n+1)x^2+(5m+1)x+1, \]
where $m,n$ are some integers. In particular, if $X$ is not the unit matrix and the order of $X$ is finite, 
then the order of $X$ is $5$ and the eigenvalues of X are $\exp(2\pi i/5),$ $\exp(4\pi i/5),$ $\ \exp(6\pi i/5), \ \exp(8\pi i/5)$.
\end{Cor}

\begin{proof}
We shall prove the first part. Let $\lambda_1,\lambda_2,\lambda_3,\lambda_4$ be the eigenvalues of $X$. Then the the eigenpolynomial $p(X)$ of $X$ is
\[ x^4-\left(\underset{1\le i\le j\le k\le4}{\sum}\lambda_i\lambda_j\lambda_k\right)x^3+\left(\underset{1\le i\le j\le4}{\sum}\lambda_i\lambda_j\right)x^2
-\left(\underset{1\le i\le4}{\sum}\lambda_i\right)x+1. \]
On the other hand, the the eigenpolynomial $p(X^{-1})$ of $X^{-1}$ is
\begin{align}
&x^4-\left(\underset{1\le i\le j\le k\le4}{\sum}\frac{1}{\lambda_i}\frac{1}{\lambda_j}\frac{1}{\lambda_k}\right)x^3+\left(\underset{1\le i\le j\le4}{\sum}
\frac{1}{\lambda_i}\frac{1}{\lambda_j}\right)x^2-\left(\underset{1\le i\le4}{\sum}\frac{1}{\lambda_i}\right)x+1 \nonumber \\
=&x^4-\left(\underset{1\le i\le4}{\sum}\lambda_i\right)x^3+\left(\underset{1\le i\le j\le4}{\sum}\lambda_i\lambda_j\right)x^2 \nonumber
-\left(\underset{1\le i\le j\le k\le4}{\sum}\lambda_i\lambda_j\lambda_k\right)x+1.
\end{align}
Since $X\in\text{Sp}(4,\mathbb{Z})$, $p(X)=p(X^{-1})$. So we can express $p(X)$ by $x^4+ax^3+bx^2+ax+1,$ where $a,b\in\mathbb{Z}$. 
It follows from the theorem that $a\equiv-4, \ b\equiv6 \ \text{mod} \ 5$. Hence the claim of the first part follows.

Next we shall prove the latter part. Let $\lambda$ be an eigenvalue of X. By the property for eigenvalues of elements of the symplectic group, 
$\bar{\lambda}, \ 1/\lambda, \ 1/\bar{\lambda}$ are also eigenvalues of X. 
If $1$ or $-1$ is an eigenvalue of X, its multiplicity is even. Since the order of $X$ is finite, we can express eigenvalues of X by 
$\exp(i\theta_1), \ \exp(-i\theta_1), \ \exp(i\theta_2), \ \exp(-i\theta_2) \ (0\le\theta_1,\theta_2\le\pi)$.
Then the eigenpolynomial of $X$ is
\[ x^4-2(\cos\theta_1+\cos\theta_2)x^3+2(\cos(\theta_1+\theta_2)+\cos(\theta_1-\theta_2)+1)x^2-2(\cos\theta_1+\cos\theta_2)x+1. \]
By the claim of the first part of the Corollary, we have
\[ -2(\cos\theta_1+\cos\theta_2)=5m+1, \ 2(\cos(\theta_1+\theta_2)+\cos(\theta_1-\theta_2)+1)=5n+1, \ m,n \in \mathbb{Z}. \]
By the addition theorem, we have
\[ 2(\cos\theta_1+\cos\theta_2)=-5m-1, \ 4\cos\theta_1\cos\theta_2=5n-1. \]
It follows from $-4\le2(\cos\theta_1+\cos\theta_2)\le4$ that $m=0,-1.$ If $m=-1$, then $\cos\theta_1,\cos\theta_2=1$ and all eigenvalues of $X$ are $1$. 
Since the order of $X$ is finite, $X$ is the unit matrix. It contradicts the assumption that $X$ is not the unit matrix. Hence $m=0$ and 
\[ \cos\theta_1+\cos\theta_2=-1/2. \]
It follows from from $-4\le4\cos\theta_1\cos\theta_2\le4$ that $n=0,1$. If $n=1$, then $\cos\theta_1=\pm1$, $\cos\theta_2=\pm1$. 
It contradicts the fact that $\cos\theta_1+\cos\theta_2=-1/2$. Hence $n=0$ and 
\[ \cos\theta_1\cos\theta_2=-1/4. \]
Combining these two equations, we have
\[ \cos^2\theta_1+1/2\cos\theta_1-1/4=0. \]
When we solve this equation for $\cos\theta_1$, 
\[ \cos\theta_1=\frac{-1\pm\sqrt{5}}{4}, \ \sin\theta_1=\frac{\sqrt{10\pm2\sqrt{5}}}{4}, \]
\[ \cos\theta_2=\frac{-1\mp\sqrt{5}}{4}, \ \sin\theta_2=\frac{\sqrt{10\mp2\sqrt{5}}}{4}. \]
Then we can verify easily that $(\exp(i\theta_1))^5, \ (\exp(i\theta_2))^5=1$. 
Hence $(\theta_1,\theta_2)=(2\pi/5,4\pi/5)$ or $(4\pi/5,2\pi/5).$
\end{proof}

\section{A relation to the other result}

In this section, we shall compare the main result of this paper with the result of Chen, Yang and Yui. 
In [2], they find the congruence subgroup $\Gamma(5,5)$ which contains the global monodromy $\Gamma$. 
Combining their result and our theorem, we can find a smaller group which contains $\Gamma$.

The congruence subgroup $\Gamma(5,5)$ is defined by 
\[ \Gamma(5,5)=\left\{ X\in\text{Sp}(4,\mathbb{Z}) \ \left| \ \gamma\equiv \begin{pmatrix}
																				1 & * & * & * \\
																				0 & 1 & * & * \\
																				0 & 0 & 1 & 0 \\
																				0 & 0 & * & 1
																					\end{pmatrix} \ \ \text{mod} \ 5 \right\}.\right. \]
$\Gamma(5,5)$ contains the principal congruence group $\Gamma(5)=\text{Ker}(\text{Sp}(4,\mathbb{Z})\to\text{Sp}(4,\mathbb{Z}/5\mathbb{Z}))$ as 
a normal subgroup of finite index. 

Let $X\in\Gamma(5,5)$ and express $X$ by
\[ \begin{pmatrix}
			5x_{11}+1&x_{12}&x_{13}&x_{14} \\
			5x_{21}&5x_{22}+1&x_{23}&x_{24} \\
			5x_{31}&5x_{32}&5x_{33}+1&5x_{34} \\
			5x_{41}&5x_{42}&x_{43}&5x_{44}+1 
				\end{pmatrix}, \ x_{ij}\in\mathbb{Z} \ (1\le i,j \le 4). \]
Then we have
\[ \text{GL}(4,\mathbb{Z})\ni P^{-1}XP\equiv\begin{pmatrix}
				1&-9x_{31}&-x_{12}+3x_{32}&-x_{14}+3x_{34} \\
				0&1&-2x_{12}&-2x_{14} \\
				0&0&1&x_{24} \\
				0&0&0&1 
					\end{pmatrix} \ \ (\text{mod} \ 5). \]
By the main theorem, if $X\in\Gamma$, then $\rho(P^{-1}XP)\in\tilde{\Gamma}$ and
\[ -9x_{31}\equiv n, \ -2x_{12}\equiv n, \ -x_{12}+3x_{32}\equiv 3n^2+2n \ \ (\text{mod} \ 5). \]
where $n$ is some integer. From a simple calculation, the above equation is equivalent to 
\[ x_{31}\equiv 3x_{12}, \ x_{32}\equiv 4x_{12}^2+4x_{12} \ \ (\text{mod} \ 5). \]
So, let $\tilde{\Gamma}(5,5)$ be
\begin{align*}
&\left\{ \left. \begin{pmatrix}
			5x_{11}+1&x_{12}&x_{13}&x_{14} \\
			5x_{21}&5x_{22}+1&x_{23}&x_{24} \\
			5x_{31}&5x_{32}&5x_{33}+1&5x_{34} \\
			5x_{41}&5x_{42}&x_{43}&5x_{44}+1 
				\end{pmatrix}\in\text{Sp}(4,\mathbb{Z}) \ \right| \ \begin{array}{l}
																		x_{31}\equiv 3x_{12}, \\
																		x_{32}\equiv 4x_{12}^2+4x_{12} \\
																		(\text{mod} \ 5)
																			\end{array} \right\}.
\end{align*}
Then we have the following Corollary.

\begin{Cor} \ 

\textup{(i)} $\tilde{\Gamma}(5,5)$ is a subgroup of $\Gamma(5,5)$.

\textup{(ii)} $\Gamma\subset\tilde{\Gamma}(5,5)\varsubsetneq\Gamma(5,5)$.

\textup{(iii)} $\tilde{\Gamma}(5,5)$ is a congruence subgroup of $\textup{Sp}(4,\mathbb{Z})$ of finite index.

\end{Cor}

\begin{proof}
Let $\rho':\Gamma(5,5)\to\text{GL}(4,\mathbb{Z}), \ X\mapsto P^{-1}XP$ and let 
$\pi=\rho\circ\rho':\Gamma(5,5)\to\text{GL}(4,\mathbb{Z}/5\mathbb{Z})$. $\tilde{\Gamma}(5,5)=\pi^{-1}(\tilde{\Gamma})$ follows from what we just mentioned. 
Since $\pi$ is a group homomorphism, $\pi^{-1}(\tilde{\Gamma})$ is a subgroup of $\Gamma(5,5)$. Hence the claim of (i) follows.

We can verify easily that $A$ and $T$ belong to $\tilde{\Gamma}(5,5)$. Therefore $\tilde{\Gamma}(5,5)$ contains $\Gamma$. 
We shall show that $\tilde{\Gamma}(5,5)$ is a proper subgroup of $\Gamma(5,5)$. 
We take $X=\begin{pmatrix}
																						1&0&0&0 \\
																						0&1&0&0 \\
																						5&0&1&0 \\
																						0&0&0&1
																							\end{pmatrix}$.
Then $X\in\Gamma(5)\subset\Gamma(5,5)$ and $X\notin\tilde{\Gamma}(5,5)$. Hence the claim of (ii) follows.

Finally, we shall show the claim of (iii). $\tilde{\Gamma}(5,5)$ contains the principal congruence subgroup 
$\Gamma(25)=\text{Ker}(\text{Sp}(4,\mathbb{Z})\to\text{Sp}(4,\mathbb{Z}/25\mathbb{Z}))$ as a normal subgroup. Hence we obtain 
$|\;\tilde{\Gamma}(5,5):\text{Sp}(4,\mathbb{Z})\;|<|\;\Gamma(25):\text{Sp}(4,\mathbb{Z})\;|=|\;\text{Sp}(4,\mathbb{Z}/25\mathbb{Z}))\;|<\infty$.
\end{proof}

\bigskip

\begin{acknowledgement}
The author would like to thank Professors Sampei Usui, Atsushi Takahashi and Keiji Oguiso for their helpful advices and suggestions. 
\end{acknowledgement}

\bigskip

\bigskip

\bigskip

\address{
Department of Mathematics \\ 
Graduate School of Science \\
Osaka University \\
Osaka 560-0043 \\
Japan
}
{ \\
u785251a@ecs.cmc.osaka-u.ac.jp}

\end{document}